\documentclass [12pt]{article}
\usepackage{amssymb,amsmath,amsthm}
\usepackage{comment}

\def\N{\mathbb{N}}

\def\F{\mathbb{F}}

\newtheorem{theorem}{Theorem}[section]
\newtheorem{proposition}[theorem]{Proposition}
\newtheorem{corollary}[theorem]{Corollary}
\newtheorem{lemma}[theorem]{Lemma}

\newtheorem{remark}[theorem]{Remark}

\begin{document}
\title{All unitary perfect polynomials over $\F_2$ with less than five distinct prime factors}
 \author{Luis H. Gallardo - Olivier Rahavandrainy \\
Department of Mathematics, University of Brest,\\
6, Avenue Le Gorgeu, C.S. 93837, 29238 Brest Cedex 3, France.\\
 e-mail : luisgall@univ-brest.fr - rahavand@univ-brest.fr}
\maketitle
\begin{itemize}
\item[a)]
Running head:  binary unitary perfect polynomials.
\item[b)]
Keywords: Sum of divisors, unitary divisors,
polynomials, finite fields,\\ characteristic $2.$
\item[c)]
Mathematics Subject Classification (2000): 11T55, 11T06.
\item[d)]
Corresponding author: Luis H. Gallardo.
\end{itemize}

\newpage

{\bf{Abstract}} We find all unitary perfect polynomials over the
prime field $\F_2$ with less than five distinct prime factors.

\section{Introduction}
Let $p$ be a prime number and let $\F_q$ be a finite field of
characteristic $p$ and order $q.$ Let $A \in \F_q[x]$ be a monic
polynomial. We say that a divisor $d$ of $A$ is unitary if $d$ is
monic and $\displaystyle{\gcd(d,\frac{A}{d}) = 1}$. Let $\omega(A)$
denote the number of distinct monic irreducible factors of $A$ over
$\F_q$ and let $\sigma(A)$ (resp. $\sigma^*(A)$) denote the sum of
all monic divisors (resp. unitary divisors) of $A$ ($\sigma$ and
$\sigma^*$ are multiplicative functions).\\
 The analogue notion over the positive integers
is the notion of unitary perfect numbers.  Only few results are known for them
(see \cite{Goto, Graham, Wall}), namely, all are even numbers, we know only five of them.
Graham  \cite{Graham} characterized three of them, namely
$6,60,87360.$ Goto \cite{Goto} proved an explicit exponential upper bound
in $k=\omega(n)$  for $n$ unitary perfect. Wall \cite{Wall} improved a previous result of Subbarao,
by proving that $\omega(n) \geq 9$
for any unitary perfect number $n.$
 
We call \emph{even} a polynomial $A$
with some zero in $\F_q,$ and \emph{odd} a polynomial that is not even.
We assume that $A \notin \F_q.$

Since $A$ and $\sigma(A)$ have the same degree it follows that
$A$ divides $\sigma(A)$ is equivalent to $\sigma(A)=A.$
If $\sigma(A) = A$ (resp. $\sigma^*(A) = A$), then we say that $A$
is a perfect (resp. unitary perfect) polynomial. We may consider
the perfect polynomials as a polynomial analogue of the multiperfect numbers.
E. F. Canaday, the
first doctoral student of Leonard Carlitz, began in $1941$ \cite{Canaday} the study
of perfect polynomials by working on the prime field $\F_2.$
 Later, in the seventies, J. T. B. Beard Jr. et al.
extended this work in several directions (see e.g. \cite{BeardU},
\cite{Beard2}, \cite{Beard}) including the study of unitary perfect polynomials.

We became interested in this subject a few years ago and obtain some
results (\cite{Gall-Rahav},
\cite{Gall-Rahav3}, \cite{Gall-Rahav4}, \cite{Gall-Rahav2},
\cite{Gall-Rahav7}, \cite{Gall-Rahav5}, \cite{Gall-Rahav6} and
\cite{Gall-Rahav8}) including  for $q \in \{2,4\}$ a complete classification of the perfect polynomials $A$
for which $\omega(A)$ is small.

We began the study of unitary perfect polynomials by considering the splitting case
when $q=p^2$ (see \cite{Gall-Rahav9}).
In this paper we study more general unitary perfect
polynomials $A$ improving on previous results of Beard et al. \cite{BeardU2} and
Beard \cite{BeardU}.   In particular we prove that $A$ must be even, contrary to perfect polynomials
for which we do not know whether or not there exist odd perfect polynomials.
More precisely, we determine here
all unitary perfect polynomials $A$, over $\F_2$, such that
$\omega(A) \leq 4$. As usual $\N$ denotes the nonnegative integers and $\N\sp{*}$ the positive
integers.
\\

Our main results are the following:\\

Let $A$ be a nonconstant polynomial over $\F_2$ such that $\omega(A)
\leq 4$, then $A$ is unitary perfect if and only if either $A$ or
$A(x+1)$ is of the form $B^{2^n}$ for some $n \in \N$ where:
$$\begin{array}{l}
- \text{ if } \omega(A) \leq 3: \\
\\
B = x(x+1),\\
B = x^3(x+1)^3(x^2+x+1)^2,\\
B(x) \in \{x^3(x+1)^2(x^2+x+1), \ x^5(x+1)^4(x^4+\cdots+x+1)\}\\
\\
- \text{ if $\omega(A)=4$:} \\
\\
i) \ B = x^6 (x+1)^4 (1+x+x^2)^3(1+x+x^4),\\
ii) \ B = x^{13} (x+1)^{8} (1+x+x^2)^4(1+x+\cdots+x^{12}),\\
iii) \ B = x^{11} (x+1)^{8} (1+x+\cdots+x^4)^2(1+x+\cdots+x^{10}),\\
iv) \ B = x^9 (x+1)^4 (1+x+x^2)^2(1+x^3+x^6),\\
v) \ B = x^{25} (x+1)^{16} (1+x+\cdots+x^4)^4(1+x^5+x^{10}+ x^{15}+
x^{20}),\\
vi) \ B = x^7 (x+1)^4 (1+x^2+x^3)(1+x+x^3),\\
vii) \ B = x^3 (x+1)^3 (1+x+x^2)^3(1+x+x^4),\\
viii) \ B = x^5 (x+1)^6 (1+x+x^2)^2(1+x +\cdots+x^4),\\
ix) \ B = x^5 (x+1)^5 (1+x^3+x^4)(1+x +\cdots+x^4),\\
x) \ B = x^{13} (x+1)^{12}(1+x+x^2)^8(1+x+\cdots+x^{12}),\\
xi) \ B = x^9 (x+1)^6 (1+x+x^2)^4(1+x^3+x^6),\\
xii) \ B = x^7 (x+1)^7 (1+x+x^3)^2(1+x^2+x^3)^2.
\end{array}$$
We may consider the family $\{x^{2^n}(x+1)^{2^n}: n \in \N\}$ as an
analogue of the family $\{x^{2^n+1}(x+1)^{2^n+1}\}$ of  trivial even perfect polynomials over $\F_2.$ \\
Note that Beard \cite{BeardU} and Beard et al. \cite{BeardU2}
computed the above list with the exception of v), x), and xi) that
are new.

Moreover, compared to the list of all perfect polynomials $A$ over $\F_2$ with $\omega(A) <5$
given in \cite{Gall-Rahav5}, we obtain an additional family of
irreducible divisors of unitary perfect polynomials:
$$\begin{array}{l} S_1(x) = 1+x^3+x^6, \ S_1(x+1),\\
S_2(x) = 1 + x^5 + x^{10} +
x^{15} + x^{20}, \ S_2(x+1)\\
S_3(x) = 1+x+\cdots+ x^{10}, \ S_3(x+1),\\
S_4(x) =  1+x+\cdots+ x^{12}, \ S_4(x+1).
\end{array}$$

It is clear from the above results that the classification of all perfect or unitary perfect
polynomials $A$ with a moderately large number $\omega(A)$  of distinct prime factors
may become very complicated.  New tools need to be discovered to make more progress
in this area.

\section{Preliminary}\label{preliminaire0}
We need the following results. Some of them are obvious, so we omit
to give their proofs. Our first result give information on the sizes of the primary
parts of unitary perfect polynomials.

\begin{lemma} \emph{(see also \cite[Theorem 1]{BeardU})}
\label{nombreminimal} If $A = \displaystyle{P_1^{h_1} \cdots
P_r^{h_r} Q_1^{k_1} \cdots Q_s^{k_s}}$ is a nonconstant unitary
perfect polynomial over $\F_q$ such that: $$\left\{\begin{array}{l}
P_1,\ldots, P_r, Q_1, \ldots, Q_s \text{are both irreducible}\\
h_1\deg(P_1)= \cdots=h_r \deg(P_r) < k_1 \deg(Q_1) \leq \cdots \leq
k_s \deg(Q_s). \end{array}\right.$$ Then:
$$r \equiv 0 \ ({\rm{mod}} p).$$
\end{lemma}

\begin{proof}
By definition, one has: $\displaystyle{0 =
\sigma^*(A)-A = \frac{A}{P_1^{h_1}} + \cdots + \frac{A}{P_r^{h_r}}+
\cdots}$ \\
In particular, $r = 1+\cdots+1$, which is the leading coefficient of
$\displaystyle{\frac{A}{P_1^{h_1}} + \cdots + \frac{A}{P_r^{h_r}}}$,
 equals $0$ in $\F_p$.
\end{proof}

\begin{lemma} \label{multiplicativity}
If $A = A_1A_2$ is unitary perfect over $\F_2$ and if $\gcd(A_1,A_2)
= 1$. Then $A_1$ is unitary perfect if and only if $A_2$ is unitary
perfect.
\end{lemma}
\begin{lemma} \label{translation}
If $A(x)$ is unitary perfect over $\F_2$, then the polynomials
$A(x+1)$ and $A^{2^n}$ are also unitary perfect over~$\F_2$, for any
$n \in \N$.
\end{lemma}

We recall here some useful notation and results in Canaday's paper \cite{Canaday}:

\begin{itemize}
\item
We define as the inverse of a polynomial $P(x)$ of degree $m$, the
polynomial $\displaystyle{P^*(x) = x^m P(\frac{1}{x})}$.
\item
We say that $P$ inverts into itself if $P = P^*$.
\item
A polynomial $P$ is complete if $P = 1 + x + \cdots + x^h$, for
some $h\in \N$.
\end{itemize}

Part iii) of the  following lemma is essentially a result of Dickson (see \cite[Lemma 2]{Canaday})

\begin{lemma} [\rm{see \cite[lemma 7]{Canaday}, \cite[Lemma 2.1]{Gall-Rahav5}}] \label{complete1}
\emph{i)} Any complete polynomial inverts into itself. \emph{ii)} If $1 + x + \cdots + x^h = PQ$, where
 $P, Q$ are irreducible,
then either $(P = P^*, Q = Q^*)$ or $(P = Q^*, Q = P^*)$.\\
\emph{iii)} If $P = P^*$, $P$ irreducible and if $P = x^a(x+1)^b+1$, then:
$$P \in \{1+x+x^2, 1 + x + \cdots + x^4\}.$$
\end{lemma}

\begin{lemma} \emph{(see \cite[Lemmata 4, 5, 6 and Theorem 8]{Canaday})}
 \label{complete2}
Let $P, Q \in \F_2[x]$ such that $P$ is irreducible and let $n,m \in \N$.\\
\emph{i)}If $1 + P + \cdots + P^{2n} = Q^m$, then $m
\in \{0,1\}$.\\
\emph{ii)} If $1 + P + \cdots + P^{2n} = Q^m A$, with $m > 1$ and $A \in
\F_2[x]$ is nonconstant, then ${\rm{deg}}(P)
> {\rm{deg}}(Q)$.\\
\emph{iii)} If $1 + x + \cdots + x^{2n} = PQ$ and $P = 1 + (x+1) + \cdots +
(x+1)^{2m}$, then $n=4$, $P= 1+x+x^2$ and $Q = P(x^3) = 1+x^3+x^6$.\\
\emph{iv)} If any irreducible factor of $1+x+\cdots + x^{2n}$ is of the
form $x^a(x+1)^b+1$, then $n \in \{1,2,3\}$.\\
\emph{v)} If $1 + x + \cdots + x^{h} = 1 + (x+1) + \cdots + (x+1)^{h}$,
then $h= 2^n-2$, for some $n \in \N$.
\end{lemma}

\begin{lemma} \label{complete3}
If $1+x+x^2$ divides $1+x+\cdots + x^{h}$, then $h \equiv 2 \mod 3$.\\
If $1+x+\cdots+x^4$ divides $1+x+\cdots + x^{h}$, then $h \equiv 4
\mod 5$.
\end{lemma}

As a special case of \cite[Theorem 2.47]{Rudolf}, we have
\begin{lemma} \label{irreduc0}
The polynomial $1+x+\cdots+x^m$  is irreducible over $\F_2$ if and
only if: $$\text{$m+1$ is a prime number and $2$ is a primitive root
in $\F_{m+1}$.}$$
\end{lemma}

Consequently one gets

\begin{lemma} \label{irreduc}
\emph{i)} The polynomial $Q(x)=1+x^5+\cdots+(x^5)^l$ is irreducible over
$\F_2$
if and only if $l = 4$.\\
\emph{ii)} The polynomial $Q(x)=1+x+\cdots+x^{3.2^r}$ is irreducible over
$\F_2$
if and only if $r = 2$.\\
\emph{iii)} The polynomial $Q(x)=1+x+\cdots+x^{5.2^r}$ is irreducible over
$\F_2$ if and only if $r = 1$.
\end{lemma}

\begin{proof}
We prove only necessity. Sufficiency is obtained by direct computations.\\
i): For $k \in \N^*$, let $\Phi_k$ be the $k$-th cyclotomic
polynomial over $\F_2$. Recall that if $k$ is a prime number, then
$\Phi_k(x) = 1+x+\cdots+x^{k-1}$.\\
If $Q(x)$ is irreducible, then $1+x+\cdots+x^l$ is also
irreducible.\\
Thus, by Lemma \ref{irreduc0}, $l+1$ is a prime number and $Q(x) = \Phi_{l+1}(x^5)$.\\
It remains to observe that if $5 \not= l+1$, then: $$\Phi_{l+1}(x^5)
= \Phi_{l+1}(x) \ \Phi_{5(l+1)}(x).$$ So that $Q$ is not irreducible
in that case. We conclude that $l=4$.\\
ii): If $Q(x)$ is irreducible, then by Lemma \ref{irreduc0}, $p=3 \
. \ 2^r + 1$ is a prime number and $2$ is a primitive root in
$\F_{p}$. So, $2$ is not a square in $\F_{p}$. By considering the
Legendre Symbol $\displaystyle{(\frac{2}{p}) =
(-1)^{\frac{p^2-1}{8}}}$, we see that we must have $r \in
\{1,2\}$.\\
The case $r=1$ does not happen since $Q(x)$ is irreducible.\\
iii): As above, we obtain: $r \in \{1,2\}$. The case $r =2$ does not
happen since $5.2^r + 1$ is prime.
\end{proof}

We prove now  the non-existence of odd unitary perfect polynomials:

\begin{lemma} \label{noodd}
Any nonconstant unitary perfect polynomial over $\F_2$ is divisible by $x$ and by
$x+1$. In particular, there is no odd unitary perfect polynomial over
$\F_2$.
\end{lemma}

\begin{proof}
If $P$ is an odd prime polynomial over $\F_2$, then
$P(0) = P(1) = 1$, so that for any positive integer $h$, $1+P(0)^h =
1+P(1)^h = 0$. Thus, the monomials $x$ and $x+1$ divide $1+P^h$.
Now, let $A$ be an unitary perfect polynomial. We have $\omega(A) \geq 2$.
If both $x, x+1$ divide $A$, then we are done.
If there exists an odd polynomial  $P \in \F_2[x]$ such that $P^h \ | \
A$ and $P^{h+1} \nmid A$, then $\sigma^*(P^h) = 1+P^h$ divides
$\sigma^*(A) = A$. So $x, x+1$ divide $A$.
\end{proof}

\begin{remark}
\begin{itemize}
\item
In the rest of the paper, we put $\overline{S}(x) = S(x+1)$ for
$S \in \F_2[x]$.
\item
For Theorems \ref{casomegainf3} and \ref{casomega4}, we shall
prove only necessity, since sufficiency is always obtained by direct
computations.
\end{itemize}
\end{remark}

\section{Case $\omega(A) \leq 3$} We prove the
following result:
\begin{theorem} \label{casomegainf3}
Let $A \in \F_2[x]$ be a polynomial such that $\omega(A) \leq 3$,
then $A$ is unitary perfect over $\F_2$ if~and~only~if either $A$ or
$\overline{A}$ is of the form $B^{2^n} \text{ for some } n \in \N,
\text {where:}$
$$\left\{\begin{array}{l}
i) \ B = x^2+x,\\
ii) \ B \in \{x^3(x+1)^2(x^2+x+1), \
x^5(x+1)^4(x^4+\cdots+x+1)\},\\
iii) \ B = x^3(x+1)^3(x^2+x+1)^2.
\end{array}
\right.$$
\end{theorem}

\subsection{Case $\omega(A) = 2$}
The following proposition gives the first part of Theorem
\ref{casomegainf3}.

\begin{proposition} Let $A \in \F_2[x]$ such that $\omega(A) =
2$, then $A$ is unitary perfect over $\F_2$ if and only if $A$ is of
the form $(x^2+x)^{2^n}$, for some $n \in \N$.
\end{proposition}

\begin{proof}
It remains to prove necessity since sufficiency is obvious.\\
The case where $A \in \{ x^hP^k, (x+1)^h P^k\}$, with $P$ odd, is
impossible by Lemma \ref{noodd}. So $A$ splits: $A =x^h(x+1)^k$. We
must have: $1+x^h = (x+1)^h, \ 1+(x+1)^k = x^k$. Hence, $h=k=2^n$,
for some $n\in \N$.
\end{proof}

Consequently the unitary perfect polynomials $A$ with $\omega(A)=2$
are exactly the perfect polynomùials with $\omega(A)=2.$

\subsection{Case $\omega(A) = 3$} In this case, $A$ is of the form
$x^{h_1}(x+1)^{k_1} P^l$, with $P$ odd.

\begin{lemma} \label{exponentodd}
If $A = x^{h_1}(x+1)^{k_1} P^l$ is an unitary perfect polynomial over
$\F_2$, then $l = 2^n$, for some nonnegative integer $n$.
\end{lemma}

\begin{proof}
Put: $l = 2^n u$, where $u$ is odd and $n \in \N$. Since the only
prime divisors of $A = \sigma^*(A)$ are $x,x+1$ and $P$, and since
$P$ does not divide $1+P^l$, the polynomial $1+P^l = \sigma^*(P^l)$
must be of the form $x^a(x+1)^b$. Thus,
$$(1+P)(1+P+ \cdots + P^{u-1}) = 1+P^u =x^c(x+1)^d.$$ Since $x, x+1$
divide $1+P$ and since $\gcd(1+P, 1+P+ \cdots + P^{u-1}) = 1$, we
conclude that $u-1 = 0$.
\end{proof}

Put $h_1 =
2^h c, \ k_1 = 2^k d$ with $c,d$ odd.
Since $A$ is unitary perfect, we have
\begin{equation}
\label{starzero}
 \left\{\begin{array}{l}
1+x^{h_1} = (x+1)^{2^h}(1+x+\cdots +x^{c-1})^{2^h},\\
1+(x+1)^{k_1} = x^{2^k}(1+(x+1)+\cdots +(x+1)^{d-1})^{2^k},\\
1+P^{2^n} = (1+P)^{2^n} = (x^{a_3}(x+1)^{b_3})^{2^n}.
\end{array}
\right.
\end{equation}

 Lemma \ref{complete2}-i) implies that:
$$1+x+\cdots
+x^{c-1}, \ 1+(x+1)+\cdots +(x+1)^{d-1} \in \{1, P\}.$$ Since $h_1$
and $k_1$ play symmetric roles and since $P$ must appear in the
right hand side of \eqref{starzero}, we may reduce the study to the two cases:
$$\begin{array}{l}
{\rm{(I)}}: 1+x+\cdots +x^{c-1}=P, \ d = 1, \\
{\rm{(II)}}: 1+x+\cdots +x^{c-1} = P = 1+(x+1)+\cdots +(x+1)^{d-1}.
\end{array}$$

\subsubsection{Case (I)}
According to Lemma \ref{complete1}-iii), we have: $P \in \{1+x+x^2,
1+x+\cdots+x^4\}$ and $c \in \{3,5\}$.\\
By considering exponents and degrees, System \eqref{starzero} implies
$$\begin{array}{l}
k=h+1, n=h \ \text{ if } c = 3,\\
k=h+2, n=h \ \text{ if } c = 5.
\end{array}$$
We obtain part ii) of Theorem \ref{casomegainf3}.

\subsubsection{Case (II)}
We have $c=d$ and $P = \overline{P}$. So, by Lemma \ref{complete1},
$P = 1+x+x^2$, and hence $c=d=3$. System  \eqref{starzero} implies:
$k=h, \ n=h+1$, and we obtain part iii) of
Theorem \ref{casomegainf3}. This completes the proof of Theorem \ref{casomegainf3}.\\

It turns out that we can also get Theorem \ref{casomegainf3}.
 as a consequence of a nice result of Swan:

\subsubsection{Another proof using Swan's Lemma}
We would like to give, here, another proof of parts ii) and iii) of
Theorem~\ref{casomegainf3}, by using Lemma \ref{nombreminimal} and
the following result about reducibility of a binary polynomial in
$\F_2[x]$:
\begin{lemma} [see \cite{Swan}, p. 1103, line 3] \label{swan}
Let $n, k \in \N$ be such that $8n > k$, then the polynomial $x^{8n} +
x^k +1$ is reducible over $\F_2$.
\end{lemma}

From that, we obviously obtain the

\begin{corollary} \label{corollaryswan}
Let $r$ be a positive integer, then the polynomial
$$P=x^{2^r} + x^{2^r-1}+1$$ is irreducible over $\F_2$ if and only if $r
\in\{1,2\}$.
\end{corollary}

We recall that $A$ is of the form $x^{h_1}(x+1)^{k_1} P^l$, with $P$
odd and $l = 2^n$ for some $n \in \N$. Put $p = \deg(P)$.
By Lemma \ref{nombreminimal}, we have either $({h_1} = {k_1} \leq
lp)$ or $({h_1}=lp \leq {k_1})$ or $({k_1}=lp \leq {h_1})$. The
third case is similar to the second
since ${h_1}$ and ${k_1}$ play symmetric roles. \\
\\
\underline{Case ${h_1} = {k_1} \leq lp$}\\
\\
We obtain $A = x^{h_1}(x+1)^{h_1}P^{2^n}$, $h_1 \leq 2^n p$. Since
$A$ is unitary perfect, we have
$$\begin{array}{l}
1+x^{h_1} = (x+1)^{b_1}P^{c_1},\\
1+(x+1)^{h_1} = x^{a_2}P^{c_2},\\
1+P^{2^n} = (1+P)^{2^n} = (x^{a_3}(x+1)^{b_3})^{2^n}.
\end{array}$$
Hence: $$\begin{array}{l}
P = x^{a_3}(x+1)^{b_3} + 1,\\
(x+1)^{b_1}P^{c_1} = 1+x^{h_1} = 1+(x+1+1)^{h_1} =
(x+1)^{a_2}{(P(x+1))}^{c_2}. \end{array}$$ It follows that:
$$a_2 = b_1, \ c_2 = c_1 \geq 1, \ P(x) = P(x+1).$$ Thus,
$c_2 = c_1=2^{n-1}$ and $a_3 = b_3$. The irreducibility of $P$
implies $a_3=b_3=1$. So, $P = x^2+x+1$. Put ${h_1} = 2^h c$, where
$c$ is odd. We have now:
$$(1+x)^{2^h}(1+x+\cdots+ x^{c-1})^{2^h} = 1+x^{h_1} =
(x+1)^{b_1}(x^2+x+1)^{2^{n-1}}.$$ Thus $c=3$ and $h=n-1$.
We get $A = B^{2^{n-1}}$, where $B=x^3(x+1)^3(x^2+x+1)^2$. So we
obtain part iii) of
Theorem \ref{casomegainf3}.\\
\\
\underline{Case ${h_1}=lp \leq {k_1}$}\\
\\
We obtain now: $A = x^{h_1}(x+1)^{k_1}P^{2^n}$, $h_1 = 2^n p \leq k_1$.
Since $A$ is unitary perfect, we have
$$\begin{array}{l}
1+x^{h_1} = (1+x^p)^{2^n}= ((x+1)^{b_1}P^{c_1})^{2^n},\\
1+(x+1)^{k_1} = x^{a_2}P^{c_2},\\
1+P^{2^n} = (1+P)^{2^n} = (x^{a_3}(x+1)^{b_3})^{2^n}.
\end{array}$$
Hence:
$$a_2+c_2p = k_1,\ b_1+c_1p = p, \ 2^nc_1+c_2 = 2^n.$$
It follows that $c_1 \in \{0,1\}$. If $c_1 = 0$, then $b_1 = p$ and $1+x^p = (x+1)^p$, so $p = 2^r$,
for some $r \in \N^*$. Thus, $a_3 + b_3 = 2^r$. Since $P =
x^{a_3}(x+1)^{b_3} +1$ is irreducible, $a_3$ and $b_3$ must be both
odd. Moreover, $c_2 = 2^n$ and $$a_2+2^n \ 2^r = a_2+c_2p = k_1 =
2^n(b_1+b_3) = 2^n(2^r + b_3).$$ Hence $$a_2 = 2^nb_3,$$ and
$$(1+ (x+1)^{2^r+b_3})^{2^n} = 1+(x+1)^{k_1} = x^{a_2}P^{c_2} =
(x^{b_3}P)^{2^n}.$$ It follows that:
$$1+ (x+1)^{2^r+b_3} = x^{b_3}P = x^{b_3}(x^{a_3}(x+1)^{b_3} +1).$$ Thus,
$$b_3 = 1,\ a_3 = 2^r - 1, \ k_1 = 2^n(2^r+1), \ P = x^{2^r - 1}(x+1)+1,$$ and $$A =
(x^{2^r}(x+1)^{2^r+1}P)^{2^n}.$$ So by Corollary
\ref{corollaryswan}, we get $r \in\{1,2\}$ and $\overline{A}$ satisfies part ii) of
Theorem~\ref{casomegainf3}.\\
If $c_1 = 1$, then $c_2 = b_1 = 0$. It follows that $1+x^p = P$,
with $p \geq 2$. This contradicts the fact that $P$ is irreducible.

\section{Case $\omega(A) = 4$}
\label{case4factors} We prove the following result:
\begin{theorem} \label{casomega4}
Let $A \in \F_2[x]$ be a polynomial such that $\omega(A) = 4$, then
$A$ is unitary perfect over $\F_2$ if~and~only~if either $A$ or
$\overline{A}$ is of the form $B^{2^n} \text{ for some } n \in \N,
\text {where:}$
$$\left\{\begin{array}{l}
i) \ B = x^6 (x+1)^4 (1+x+x^2)^3(1+x+x^4),\\
ii) \ B = x^{13} (x+1)^{8} (1+x+x^2)^4(1+x+\cdots+x^{12}),\\
iii) \ B = x^{11} (x+1)^{8} (1+x+\cdots+x^4)^2(1+x+\cdots+x^{10}),\\
iv) \ B = x^9 (x+1)^4 (1+x+x^2)^2(1+x^3+x^6),\\
v) \ B = x^{25} (x+1)^{16} (1+x+\cdots+x^4)^4(1+x^5+x^{10}+ x^{15}+
x^{20}),\\
vi) \ B = x^7 (x+1)^4 (1+x^2+x^3)(1+x+x^3),\\
vii) \ B = x^3 (x+1)^3 (1+x+x^2)^3(1+x+x^4),\\
viii) \ B = x^5 (x+1)^6 (1+x+x^2)^2(1+x +\cdots+x^4),\\
ix) \ B = x^5 (x+1)^5 (1+x^3+x^4)(1+x +\cdots+x^4),\\
x) \ B = x^{13} (x+1)^{12}(1+x+x^2)^8(1+x+\cdots+x^{12}),\\
xi) \ B = x^9 (x+1)^6 (1+x+x^2)^4(1+x^3+x^6),\\
xii) \ B = x^7 (x+1)^7 (1+x+x^3)^2(1+x^2+x^3)^2.
\end{array}
\right.$$
\end{theorem}

The following proposition gives more details about the form of an
unitary perfect polynomial.

\begin{proposition} \label{exponentodd2}
Every unitary perfect polynomial $A$ over $\F_2$, with $\omega(A) =
4$, is of the form $x^{h_1}(x+1)^{k_1} P^{2^l u}Q^{2^m}$, where:
$$\begin{array}{l}
\text{i) $P, Q,u$ are odd, ${\rm{deg}}(P) \leq {\rm{deg}}(Q)$},\\
\text{ii) $h_1,k_1 \in \N^*$, $l,m \in \N$ and either $(u =1)$ or
$(u=3,
\ Q = 1+P+ P^2)$,}\\
\text{iii) $P \in \{1+x+x^2, 1+x+\cdots+x^4\}$ if $P$ is complete},\\
\text{iv) ${\rm{deg}}(Q) \geq 4$ if $\ Q$ is complete}.
\end{array}$$
\end{proposition}

\begin{proof}
 First of all, $x$ and $x+1$ divide $A$ by Lemma
\ref{noodd}. So $$A = x^{h_1}(x+1)^{k_1} P^{r}Q^{s},$$ for some
$h_1,k_1,r,s \in \N^*$. Put $r = 2^l u, \ s = 2^m v$, where $u,v$
are odd and $l,m \in \N$. Consider
$$\sigma^*(Q^s) = 1+Q^s = (1+Q)^{2^m}(1+Q+\cdots + Q^{v-1})^{2^m}.$$
Since $x$ and $x+1$ divide $1+Q$, they do not divide $1+Q+\cdots +
Q^{v-1}$. Hence, $1+Q+\cdots + Q^{v-1} \in \{1,P\}$, by Lemma
\ref{complete2}-i). If $v-1 \geq 2$, then $1+Q+\cdots + Q^{v-1} =
P$. This is impossible because ${\rm{deg}}(P) \leq {\rm{deg}}(Q)$.
Thus, $v-1 = 0$ and $s = 2^m$. Now, by considering degrees, we see
that the irreducible odd polynomial $Q$ does not divide $1+P$. It
follows that $(1+P)^{2^l}(1+P+\cdots + P^{u-1})^{2^l} = 1+P^r =
\sigma^*(P^r)$ must be of the form $x^a(x+1)^bQ^c$. Thus, by Lemma
\ref{complete2}-i):
$$1+P+\cdots +
P^{u-1} \in \{1, Q\}.$$ We conclude that either $(u = 1)$ or
$(1+P+\cdots + P^{u-1} = Q)$.\\
If $u > 1$, then put $u=2w+1$. We get
$$1+Q^{2^m} = (1+Q)^{2^m} = \text{\Large{(}}P(1+P+\cdots+P^{u-2})\text{\Large{)}}^{2^m} =$$
$$
\text{\LARGE{(}}P(1+P)\text{\Large{(}}1+P+\cdots+P^{w-1}\text{\Large{)}}^2
\text{\LARGE{)}}^{2^m}.$$ Since $x, x+1$ and $P$ divide $1+Q$ and
since $x,x+1$ divide $1+P$, none of the irreducible divisors of $A$ does
divide $1+P+\cdots+P^{w-1}$.
Hence $w=1$, $u = 3$ and $Q = 1+P+P^2$.
Since ${\rm{deg}}(P) \leq {\rm{deg}}(Q)$, the irreducible
polynomial $Q$ does not divide $1+P$. So $P$ is always of the form
$x^a(x+1)^b + 1$. If $P$ is complete, then by parts i) and iii) of
Lemma \ref{complete1}, we have $P \in \{1+x+x^2,
1+x+\cdots+x^4\}$.
Finally, if $Q$ is complete, since $1+x+x^2$ is the only degree $2$ odd
irreducible polynomial over $\F_2$,  we must have
${\rm{deg}}(Q) \geq 4$.
\end{proof}

Put
 $$p={\rm{deg}}(P), \ q = {\rm{deg}}(Q), \ h_1 = 2^hc, \ k_1 =
2^kd, \text{ with $c, d$ odd}.$$
Since $A$ is unitary perfect and
since $Q$ does not divide $1+P$, we have:
\begin{equation}
\label{star}
\ \left\{\begin{array}{l}
1+x^{h_1} = (1+x^c)^{2^h}= (1+x)^{2^h}(1+x+\cdots +x^{c-1})^{2^h} =
(1+x)^{2^h}P^{2^hc_1}Q^{2^hd_1},\\
1+(x+1)^{k_1} = x^{2^k}(1+(1+x)+\cdots +(1+x)^{d-1})^{2^k} = x^{2^k}P^{2^kc_2}Q^{2^kd_2},\\
1+P^{2^lu} = (1+P)^{2^l}(1+P+\cdots+P^{u-1})^{2^l} = (x^{a_3}(1+x)^{b_3})^{2^l}Q^{2^ld_3},\\
1+Q^{2^m} = (1+Q)^{2^m} = (x^{a_4}(1+x)^{b_4}P^{c_4})^{2^m}.
\end{array}
\right.
\end{equation}
 By considering degrees and exponents of $x, x+1, P$ and
$Q$, \eqref{star} implies:
\begin{equation}
\label{2star}
 \ \left\{\begin{array}{l}
2^hc = 2^h(1+pc_1+qd_1) = 2^k + 2^l a_3 + 2^m a_4,\\
2^kd = 2^k(1+pc_2+qd_2) = 2^h + 2^l b_3 + 2^m b_4,\\
2^lup = 2^l(a_3 + b_3+qd_3) = (2^hc_1 + 2^kc_2+2^mc_4)p,\\
2^mq = 2^m(a_4+b_4+pc_4) = (2^hd_1 + 2^kd_2+2^ld_3)q.
\end{array}
\right.
\end{equation}

 By Lemma \ref{complete2}, $c_1,d_1,c_2,d_2,d_3 \in
\{0,1\}$ so that:
$$1+x+\cdots +x^{c-1}, \ 1+(1+x)+\cdots +(1+x)^{d-1} \in \{1, P, Q,
PQ\}.$$ Since $h_1$ and $k_1$ play symmetric roles, and since $x,
x+1, P$ and $Q$ must divide $A=\sigma^*(A)$, it is sufficient to consider  the
following ten cases:
$$\begin{array}{l}
{\rm{(I)}}: c=d=1, \\
{\rm{(II)}}: 1+x+\cdots +x^{c-1} = P, \ d = 1, \\
{\rm{(III)}}: 1+x+\cdots +x^{c-1} = Q, \ d = 1, \\
{\rm{(IV)}}: 1+x+\cdots +x^{c-1} = PQ, \ d = 1,\\
{\rm{(V)}}: 1+x+\cdots +x^{c-1} = P = 1+(x+1)+\cdots +(x+1)^{d-1},\\
{\rm{(VI)}}: 1+x+\cdots +x^{c-1} = Q, \ 1+(x+1)+\cdots +(x+1)^{d-1}
=P,\\
{\rm{(VII)}}: 1+x+\cdots +x^{c-1} = PQ, \ 1+(x+1)+\cdots
+(x+1)^{d-1}
=P,\\
{\rm{(VIII)}}: 1+x+\cdots +x^{c-1} = Q = 1+(x+1)+\cdots +(x+1)^{d-1},\\
{\rm{(IX)}}: 1+x+\cdots +x^{c-1} = PQ, \ 1+(x+1)+\cdots +(x+1)^{d-1}
=Q,\\
{\rm{(X)}}: 1+x+\cdots +x^{c-1} = PQ= 1+(x+1)+\cdots +(x+1)^{d-1}.
\end{array}$$

\subsection{Case (I)}
In this case, if $u = 1$, then since $Q$ must appear in the right
hand side of System \eqref{star}, $Q$ must divide $1+P$, which is
impossible.
So, $u = 3$ and $1+Q = P(P+1)$. Thus, System  \eqref{star} implies that
$c_4 = 1$ and $3\cdot2^l = c_4 \cdot 2^m = 2^m$ so that $3$ divides $2^m$. It
is impossible.

\subsection{Case (II)}
As above, $u = 3$ and $Q = 1+P+P^2$.
By Proposition \ref{exponentodd2}, we get
$$P \in \{1+x+x^2, 1+x+\cdots+x^4\} \text{ and } c \in \{3, 5\}.$$
If $P = 1+x+\cdots+x^4$, then:
$$Q=1+P+P^2=1+x+x^3+x^6+x^8=(1+x+x^2)(1+x^2+x^4+x^5+x^6),$$ which is
reducible.\\
So we must have: $P = 1+x+x^2$. Thus, $c=3$ and $Q = 1+x+x^4$.
System \eqref{2star} implies that:
$$l = m, \ h = m+1, \ k = m+2.$$ We obtain part i) of Theorem
\ref{casomega4}.

\subsection{Case (III)}
$P$ must divide $1+Q$ since it must appear in the right hand side~of~ \eqref{star}.\\
Put: $c-1 = 2^r s$, with $s$ odd. We get
$$x^{a_4}(1+x)^{b_4+1}P^{c_4} = (1+x)(1+Q) =
x(x+1)(1+x+\cdots+x^{c-2}).$$ Thus, $a_4 = 1$ and
$$(x+1)^{b_4+1}P^{c_4} = (1+x)(1+x+\cdots+x^{c-2}) = 1+x^{c-1} =
(1+x)^{2^r}(1+x+\cdots+x^{s-1})^{2^r}.$$ We conclude that:
$$b_4 = 2^r - 1, \ c_4 = 2^r, \ P = 1+x+\cdots+x^{s-1}.$$
By Proposition \ref{exponentodd2}, we get
$$P \in \{1+x+x^2, 1+x+\cdots+x^4\}.$$
Thus, $c \in \{3 \cdot 2^r + 1, 5 \cdot 2^r + 1\}$, and by Lemma \ref{irreduc}, $c \in \{11,13\}$.
It follows that we must have $$\begin{array}{l} u=1, \ d_3=0, \\
P = 1+x+x^2, \ Q
=
1+x+\cdots+x^{12} \ \text{ if $c = 13$},\\
P = 1+x+\cdots+x^4, \ Q = 1+x+\cdots+x^{10} \ \text{ if $c = 11$}.
\end{array}$$

System \eqref{2star} implies
$$\begin{array}{l}
m=h, \ l=h+2, \ k = h+3 \ \text{ if $c = 13$},\\
m=h, \ l=h+1, \ k = h+3 \ \text{ if $c = 11$}.
\end{array}$$ We obtain parts ii) and iii) of Theorem \ref{casomega4}.

\subsection{Case (IV)}
We get $1+x+\cdots+x^{c-1} = PQ,$ and by Lemma \ref{complete1}: $P
\in \{P^*, Q^*\}$.

\subsubsection{Case $P =
P^*$} In this case, by Lemma \ref{complete1}-iii), we have: $P \in
\{1+x+x^2,1+x+\cdots +
x^4\}$.\\
$\bullet$ If $P = 1+x+x^2$, then by Lemma \ref{complete2}-iii), the
only possibility is $$c = 9, \ Q = 1+x^3+x^6.$$ So, we must have
$$u=1.$$ System \eqref{2star} implies the following:
$$m=h, \ l=h+1, \ k = h+2.$$
We obtain then part iv) of Theorem \ref{casomega4}.\\

$\bullet$ If $P = 1+x+\cdots+x^4$,
then $1+x+\cdots+x^4$ divides
$1+x+\cdots+x^{c-1}$.\\
So, by Lemma \ref{complete3}, $c$ is divisible by $5$. Put $c=5 w$.
We get $Q = 1+x^5+x^{10}+ \cdots + (x^5)^{w-1} \not= 1+P+P^2$. Thus,
by Lemma \ref{irreduc}-i) and by Proposition \ref{exponentodd2}, we
have
$$c=5w=25, \ u=1, \ P = 1+x+\cdots+x^4, \ Q = 1+x^5+ x^{10}+x^{15}+
x^{20}.$$ System  \eqref{2star} implies
$$m=h, \ l=h+2, \ k = h+4.$$ So
we obtain part v) of Theorem \ref{casomega4}.

\subsubsection{Case $P =Q^*$}
We get $p=q$. So both $P$ and $Q$ are of
the form $x^a(x+1)^b+1$.
We conclude by Lemma \ref{complete2}-iv) that:
$$c = 7, \ P, Q \in \{1+x^2+x^3, 1+x+x^3\}.$$ It follows that $Q \not=
1+P+P^2$ and $u =1$.
System \eqref{2star} implies
$$l=m=h, \ k=h+2.$$
We obtain then part vi) of Theorem \ref{casomega4}.

\subsection{Case (V)}
In this case, by Lemma \ref{complete1}-iii), $P = 1+x+x^2$ and
$c=d=3$. Moreover, $u$ must be equal to $3$. So, $Q = 1+P+P^2 =
1+x+x^4$. System \eqref{2star} implies now:
$$l=m=k=h.$$
Consequently we obtain part vii) of Theorem \ref{casomega4}.

\subsection{Case (VI)}
In this case, $\overline{P} \in \{1+x+x^2, 1+x+\cdots+x^4\}$ by
Lemma
\ref{complete2}-iv).
So $Q \not= 1+P+P^2$ and hence $u = 1$.

\subsubsection{Case where $P$ does not divide $1+Q$}
In this case, both $P$ and $Q$ are of the form $x^a(x+1)^b+1$. By
Lemma \ref{complete2}-iv) and Proposition
\ref{exponentodd2}-iii)-iv), we have two possibilities:
$$\begin{array}{l}
\overline{P} = P =
1+x+x^2, \ Q = 1+x+\cdots+x^4,\\
\overline{P}  = 1+x+\cdots+x^4 = Q.
\end{array}$$ Thus $(c,d) \in \{(5,3), (5,5)\}$.
System \eqref{2star}  implies $$\begin{array}{l} m=h, \ l=k=h+1
\text{ if $c=5, \  d=3$},\\
l=m=k=h \text{ if $c=d=5$}. \end{array}$$ We obtain parts viii) and
ix) of Theorem \ref{casomega4}.

\subsubsection{Case where $P$ divides $1+Q$}
In this case, $P$ must divide $\displaystyle{\frac{1+Q}{x} = 1+x+\cdots+ x^{c-2}}$.
Moreover, according to System \eqref{star}, we have
$$a_4 = 1, \ 1+x+\cdots+x^{c-2} = (x+1)^{b_4}P^{c_4}.$$
Thus, if we put $c-1 = 2^r s$, with $s$ odd, we obtain
$$(1+x)^{2^r}(1+x+\cdots+x^{s-1})^{2^r} = (1+x^s)^{2^r} = 1+x^{c-1} = (x+1)^{b_4+1}P^{c_4}.$$
We conclude that: $$b_4 = 2^r - 1,$$ and by Lemma
\ref{complete2}-i):
$$P = 1+x+\cdots+x^{s-1}, \ c_4 = 2^r.$$
Hence, by Lemmata \ref{complete2}-v) and \ref{complete1}-iii)the
only possiblity that remains is
$$\overline{P} = 1+x+x^2 = P, \ s=3, \ c = 3 \cdot 2^r+1.$$
It follows that $r = 2$ by Lemma \ref{irreduc}.
System \eqref{2star} implies that: $$m=h, \ k=h+2,\ l=h+3.$$ We
obtain part x) of Theorem \ref{casomega4}.

\subsection{Case (VII)}
In this case, $P$ divides $1+x+\cdots+ x^{c-1}$. By Lemma
\ref{complete2}-iii), we get $$c=9, \ P = 1+x+x^2, \ d=3, \ Q =
1+x^3+x^6.$$ Moreover, $u=1$ since $Q  \not= 1+P+P^2$. \\
System \eqref{2star}
implies that:
$$m=h, \ k=h+1,\ l=h+2.$$ We obtain part xi) of Theorem
\ref{casomega4}.

\subsection{Case (VIII)}
In this case, by Lemma \ref{complete2}-v) and by Proposition \ref{exponentodd2}-iv), $c=d=2^w-1 \geq 5$.\\
Since $P$ must appear in the right hand side of \eqref{star}
, it must
divide $1+Q =
x(1+x+\cdots+x^{c-2})$. Hence $P$ divides $1+x+\cdots+x^{c-2}$.
Thus, $$a_4 = 1 \text{ and } (x+1)^{b_4}P^{c_4}= 1+x+\cdots+x^{c-2}
= (1+x)(1+x+\cdots+x^{2^{w-1}-2})^2.$$ We deduce that:
$$b_4 = 1, \ c_4 = 2, \ P = 1+x+\cdots+x^{2^{w-1}-2}.$$
By Proposition \ref{exponentodd2}-iii), we must have $$2^{w-1}-2 \in
\{2,4\}.$$ So $w=3$ and $Q = 1+x+\cdots+x^7=(1+x)^7$ which is not
irreducible.

\subsection{Case (IX)}
In this case, $Q$ divides $1+x+\cdots+ x^{c-1}$. By Lemma
\ref{complete2}-iii), we get $$Q = 1+x+x^2, \ P = 1+x^3+x^6.$$ This
contradicts the fact: $\deg(P) \leq \deg(Q)$.

\subsection{Case (X)}
In this case, by Lemma \ref{complete2}-v), by Proposition
\ref{exponentodd2}-iv) and by Lemma \ref{complete1}-ii), we get
$$c=d=2^w-1 \geq 5, \text{ and either $(P=P^*, Q = Q^*)$
 or $(P = Q^*)$}.$$

\subsubsection{Case where $P=P^*,\ Q = Q^*$}
We have by Lemma \ref{complete1}-iii):
$P \in \{1+x+x^2,1+x+\cdots+x^4\}$.\\
$\bullet$ If $P = 1+x+x^2$, by Lemma \ref{complete2}-iii), $Q = 1+x^3+x^6$.
Thus, $c = 9 =2^w-1$. This is impossible.\\
$\bullet$ If $P = 1+x+\cdots+x^4$, then $\overline{P}$ divides
$1+x+\cdots +x^{d-1}$.
So, by Lemma \ref{complete2}, $d-1 =8$. This is
impossible.

\subsubsection{Case where $P=Q^*$}
We have $p = q$ and both $P,Q$ are of the form $x^a(x+1)^b+1$.
By Lemma \ref{complete2}-iv), $$c=d=7 \text{ and } P,Q \in
\{1+x+x^3, 1+x^2+x^3\}.$$ Moreover $u = 1$, by Proposition
\ref{exponentodd2}-ii).
System \eqref{2star}
implies that:
$$l=m=h+1, \ k = h.$$ We obtain finally part xii) of Theorem
\ref{casomega4}.
This completes the proof of the Theorem.

\section{Acknowledgments}
We are grateful to the referee of a first version of this paper
for careful reading and for suggestions that improved the presentation of the paper.
We are including in the next section his report (but excluding
the detailed technical suggestions to authors).

\section{Report on preliminary version and conclusion}
Referee report on the paper "All unitary perfect polynomials
over F2 with less than five distinct prime factors" by Luis H.
Gallardo and Olivier Rahavandrainy.\\

The authors are studying the problem of finding all the unitary perfect poly
nomials over finite fields. The present paper contains the full classification
of all the perfect unitary polynomials over F2 and serves as a continuation
of a series of their publication devoted to the same topic. Previously the
problem was studied by E.F. Canaday, J.T.B. Beard Jr, A.T. Bulloc, M.S.
Harbin, J.R. Oconnel Jr, K.I. West. The latest publication on study of the
perfect unitary polynomials was published in 1991, and this makes papers
of the mentioned authors hardly available. Moreover, publications [2]-[4] in
the reference list is unavailable since the journal Rend. Acad. Lincei they
published in has status "no longer indexed" in database of the AMS and the
journal's webpage containing the mentioned volumes was not found. Happily
the authors are citing the papers [2]-[4] only in the history of the question.
The general idea of the proofs of the results in the paper is rather ele-
mentary. But it requires a great scope of computations and applies more
deep results on irreducibility of the polynomials. Some of these irreducibil-
ity results was proved by the authors in their previous papers. In general
the paper makes good impression by numerous tricks used by the authors to
simplify computations. The paper worth to be published in the Journal, it
presents a new research which devoted to an interesting problem. The au-
thors gives several interesting ideas, combination of which solves a problem.
I found several misprints and places where arguments of proofs are not clear.
I'd like to recommend the authors to correct misprints and clarify unclear
arguments in proofs.

\section{Conclusion}
From the (seemingly favorable ?) report above it was deduced that the preliminary version
of this paper was not
suitable for publication in the IJNT.

\def\thebibliography#1{\section*{\titrebibliographie}
\addcontentsline{toc}
{section}{\titrebibliographie}\list{[\arabic{enumi}]}{\settowidth
 \labelwidth{[
#1]}\leftmargin\labelwidth \advance\leftmargin\labelsep
\usecounter{enumi}}
\def\newblock{\hskip .11em plus .33em minus -.07em} \sloppy
\sfcode`\.=1000\relax}
\let\endthebibliography=\endlist

\def\biblio{\def\titrebibliographie{References}\thebibliography}
\let\endbiblio=\endthebibliography




\newbox\auteurbox
\newbox\titrebox
\newbox\titrelbox
\newbox\editeurbox
\newbox\anneebox
\newbox\anneelbox
\newbox\journalbox
\newbox\volumebox
\newbox\pagesbox
\newbox\diversbox
\newbox\collectionbox
\def\fabriquebox#1#2{\par\egroup
\setbox#1=\vbox\bgroup \leftskip=0pt \hsize=\maxdimen \noindent#2}
\def\bibref#1{\bibitem{#1}


\setbox0=\vbox\bgroup}
\def\auteur{\fabriquebox\auteurbox\styleauteur}
\def\titre{\fabriquebox\titrebox\styletitre}
\def\titrelivre{\fabriquebox\titrelbox\styletitrelivre}
\def\editeur{\fabriquebox\editeurbox\styleediteur}

\def\journal{\fabriquebox\journalbox\stylejournal}

\def\volume{\fabriquebox\volumebox\stylevolume}
\def\collection{\fabriquebox\collectionbox\stylecollection}
{\catcode`\- =\active\gdef\annee{\fabriquebox\anneebox\catcode`\-
=\active\def -{\hbox{\rm
\string-\string-}}\styleannee\ignorespaces}}
{\catcode`\-
=\active\gdef\anneelivre{\fabriquebox\anneelbox\catcode`\-=
\active\def-{\hbox{\rm \string-\string-}}\styleanneelivre}}
{\catcode`\-=\active\gdef\pages{\fabriquebox\pagesbox\catcode`\-
=\active\def -{\hbox{\rm\string-\string-}}\stylepages}}
{\catcode`\-
=\active\gdef\divers{\fabriquebox\diversbox\catcode`\-=\active
\def-{\hbox{\rm\string-\string-}}\rm}}
\def\ajoutref#1{\setbox0=\vbox{\unvbox#1\global\setbox1=
\lastbox}\unhbox1 \unskip\unskip\unpenalty}
\newif\ifpreviousitem
\global\previousitemfalse
\def\separateur{\ifpreviousitem {,\ }\fi}
\def\voidallboxes
{\setbox0=\box\auteurbox \setbox0=\box\titrebox
\setbox0=\box\titrelbox \setbox0=\box\editeurbox
\setbox0=\box\anneebox \setbox0=\box\anneelbox
\setbox0=\box\journalbox \setbox0=\box\volumebox
\setbox0=\box\pagesbox \setbox0=\box\diversbox
\setbox0=\box\collectionbox \setbox0=\null}
\def\fabriquelivre
{\ifdim\ht\auteurbox>0pt
\ajoutref\auteurbox\global\previousitemtrue\fi
\ifdim\ht\titrelbox>0pt
\separateur\ajoutref\titrelbox\global\previousitemtrue\fi
\ifdim\ht\collectionbox>0pt
\separateur\ajoutref\collectionbox\global\previousitemtrue\fi
\ifdim\ht\editeurbox>0pt
\separateur\ajoutref\editeurbox\global\previousitemtrue\fi
\ifdim\ht\anneelbox>0pt \separateur \ajoutref\anneelbox
\fi\global\previousitemfalse}
\def\fabriquearticle
{\ifdim\ht\auteurbox>0pt        \ajoutref\auteurbox
\global\previousitemtrue\fi \ifdim\ht\titrebox>0pt
\separateur\ajoutref\titrebox\global\previousitemtrue\fi
\ifdim\ht\titrelbox>0pt \separateur{\rm in}\
\ajoutref\titrelbox\global \previousitemtrue\fi
\ifdim\ht\journalbox>0pt \separateur
\ajoutref\journalbox\global\previousitemtrue\fi
\ifdim\ht\volumebox>0pt \ \ajoutref\volumebox\fi
\ifdim\ht\anneebox>0pt  \ {\rm(}\ajoutref\anneebox \rm)\fi
\ifdim\ht\pagesbox>0pt
\separateur\ajoutref\pagesbox\fi\global\previousitemfalse}
\def\fabriquedivers
{\ifdim\ht\auteurbox>0pt
\ajoutref\auteurbox\global\previousitemtrue\fi
\ifdim\ht\diversbox>0pt \separateur\ajoutref\diversbox\fi}
\def\endbibref
{\egroup \ifdim\ht\journalbox>0pt \fabriquearticle
\else\ifdim\ht\editeurbox>0pt \fabriquelivre
\else\ifdim\ht\diversbox>0pt \fabriquedivers \fi\fi\fi
.\voidallboxes}

\let\styleauteur=\sc
\let\styletitre=\it
\let\styletitrelivre=\sl
\let\stylejournal=\rm
\let\stylevolume=\bf
\let\styleannee=\rm
\let\stylepages=\rm
\let\stylecollection=\rm
\let\styleediteur=\rm
\let\styleanneelivre=\rm

\begin{biblio}{99}

\begin{bibref}{Beard2}
\auteur{J. T. B. Beard Jr}  \titre{Perfect polynomials Revisited}
\journal{Publ. Math. Debrecen} \volume{38/1-2} \pages 5-12 \annee
1991
\end{bibref}

\begin{bibref}{BeardU}
\auteur{J. T. B. Beard Jr} \titre{Unitary perfect polynomials over
$GF(q)$} \journal{Rend. Accad. Lincei} \volume{62} \pages 417-422
\annee 1977
\end{bibref}

\begin{bibref}{BeardU2}
\auteur{J. T. B. Beard Jr, A. T. Bullock, M. S. Harbin}
\titre{Infinitely many perfect and unitary perfect polynomials}
\journal{Rend. Accad. Lincei} \volume{63} \pages 294-303 \annee 1977
\end{bibref}

\begin{bibref}{Beard}
\auteur{J. T. B. Beard Jr, J. R. Oconnell Jr, K. I. West}
\titre{Perfect polynomials over $GF(q)$} \journal{Rend. Accad.
Lincei} \volume{62} \pages 283-291 \annee 1977
\end{bibref}

\begin{bibref}{Canaday}
\auteur{E. F. Canaday} \titre{The sum of the divisors of a
polynomial} \journal{Duke Math. J.} \volume{8} \pages 721-737 \annee
1941
\end{bibref}

\begin{bibref}{Gall-Rahav}
\auteur{L. Gallardo, O. Rahavandrainy} \titre{On perfect polynomials
over $\F_4$} \journal{Port. Math. (N.S.)} \volume{62(1)} \pages
109-122 \annee 2005
\end{bibref}

\begin{bibref}{Gall-Rahav3}
\auteur{L. Gallardo, O. Rahavandrainy} \titre{Perfect polynomials
over $\F_4$ with less than five prime factors} \journal{Port. Math.
(N.S.) } \volume{64(1)} \pages 21-38 \annee 2007
\end{bibref}

\begin{bibref}{Gall-Rahav4}
\auteur{L. H. Gallardo, O. Rahavandrainy} \titre{Odd perfect
polynomials over $\F_2$} \journal{J. Th\'eor. Nombres Bordeaux}
\volume{19} \pages 165-174 \annee 2007
\end{bibref}

\begin{bibref}{Gall-Rahav2}
\auteur{L. H. Gallardo, O. Rahavandrainy} \titre{On splitting
perfect polynomials over $\F_{p^p}$} \journal{Preprint (2007)}
\end{bibref}

\begin{bibref}{Gall-Rahav7}
\auteur{L. H. Gallardo, O. Rahavandrainy} \titre{There is no odd
perfect polynomial over $\F_2$ with four prime factors}
\journal{Port. Math. (N.S.)} \volume{66(2)} \pages 131-145 \annee
2009
\end{bibref}

\begin{bibref}{Gall-Rahav5}
\auteur{L. H. Gallardo, O. Rahavandrainy} \titre{Even perfect
polynomials over $\F_2$ with four prime factors} \journal{Intern. J.
of Pure and Applied Math.} \volume{52(2)} \pages 301-314 \annee 2009
\end{bibref}

\begin{bibref}{Gall-Rahav6}
\auteur{L. H. Gallardo, O. Rahavandrainy} \titre{On splitting
perfect polynomials over $\F_{p^2}$} \journal{Port. Math. (N.S.) }
\volume{66(3)} \pages 261-273 \annee 2009
\end{bibref}

\begin{bibref}{Gall-Rahav8}
\auteur{L. H. Gallardo, O. Rahavandrainy} \titre{All perfect
polynomials with up to four prime factors over $\F_4$}
\journal{Math. Commun.} \volume{14(1)} \pages 47-65 \annee 2009
\end{bibref}

\begin{bibref}{Gall-Rahav9}
\auteur{L. H. Gallardo, O. Rahavandrainy} \titre{On unitary
splitting perfect polynomials over $\F_{p^2}$} \journal{Preprint
(2009)}
\end{bibref}

\begin{bibref}{Goto}
\auteur{Goto, Takeshi}
\titre{Upper bounds for unitary perfect numbers and unitary harmonic numbers}
\journal{Rocky Mountain J. Math.}
\volume{37(5)}
\pages 1557-1576 \annee 2007
\end{bibref}

\begin{bibref}{Graham}
\auteur{Graham, S. W.} \titre{Unitary perfect numbers with squarefree odd part}
\journal{Fibonacci Quart.} \volume{26(4)} \pages 312-317 \annee 1989
\end{bibref}

\begin{bibref}{Rudolf}
\auteur{R. Lidl, H. Niederreiter} \titrelivre{Finite Fields,
Encyclopedia of Mathematics and its applications} \editeur{Cambridge
University Press} \anneelivre 1983 (Reprinted 1987)
\end{bibref}

\begin{bibref}{Swan}
\auteur{R. G. Swan} \titre{Factorization of polynomials over finite
fields} \journal{Pacific J. Math.} \volume{12} \pages 1099-1106
\annee 1962
\end{bibref}

\begin{bibref}{Wall}
\auteur{Wall, Charles R.} \titre{New unitary perfect numbers have at least nine odd components}
\journal{Fibonacci Quart.} \volume{26(4)} \pages 312-317
\annee 1988
\end{bibref}

\end{biblio}

\end{document}